%% file: ArXives_Version.tex
\documentclass{article}
\usepackage{amsmath,amssymb,amsthm,eucal,graphics}
\usepackage[cp866]{inputenc}
\usepackage[russian]{babel}
\usepackage{amsfonts,amssymb,eucal}
\usepackage{amsbsy}
\usepackage{latexsym}
\makeatletter
\gdef\firstpage{1}

\def\frsthdr{}
\def\firstpageone{0\thepage}
\def\firstpagetwo{00\thepage}
\def\firstpagethree{000\thepage}
\def\firstpagemark{\ifnum\firstpage <10  \firstpageone
\else\ifnum\firstpage<100 \firstpagetwo \else \ifnum\firstpage
<1000 \firstpagethree \else \firstpageone\fi\fi\fi}
\def\footline{\ifnum\thepage=\firstpage \footlineone
\else\footlinetwo\fi}
\def\footlineone{}
\def\footlinetwo{}
\def\titles{}
\def\authors{}
\def\oddhedr{\ifnum\thepage=\firstpage \firsthdr \else \odhdr \fi}
\def\firsthdr{\hspace{\fill} \sl \frsthdr \hspace{\fill}\hbox{}}
\def\odhdr{\hspace{\fill}\sl\rightmark \titles \hspace{\fill}
\rm \thepage}
\def\evnhedr{\ifnum\thepage=\firstpage \firsthdr \else \evhdr \fi}
\def\evhdr{\noindent \rm \thepage\hspace*{\fill} \sl\leftmark
\authors \hspace*{\fill}\hbox{}}
\def\ps@newpstyle{\def\@oddhead{
\hspace{-0.65em} \vbox{\oddhedr\vskip 1mm \hrule width
\textwidth}}
\def\@evenhead{
\hspace{-0.65em} \vbox{\evnhedr\vskip 1mm \hrule width
\textwidth}
}\textsc{}
\def\@oddfoot{\footline}
\def\@evenfoot{\@oddfoot}}

\def\refer
{\begin{small}
\begin{enumerate}
\leftmargin=0pt \rightmargin=0pt
\itemsep=0pt
\parsep=0pt}
\def\endref{\end{enumerate}\end{small} }
\newtheorem {Lemma}       {Леммa }
\newtheorem {Theorem}     {Теоремa }

\begin{document}
\setcounter{page}{\firstpage}
\pagestyle{newpstyle}
\Russian
\sloppy
\rm
\input alglog2017rev2.txt

\end{document}

%% file: alglog2017rev2.txt
 {\large \centerline{{\bf
 Комбинаторика двоичных слов и коразмерности   }}}

{\large \centerline{{\bf  тождеств     левонильпотентных алгебр}}}

\vskip 0.2in \centerline{М. В. Зайцев \footnote{Работа первого автора
поддержана Российским Научным Фондом, грант № 16-11-10013}}

\vskip 0.2in \centerline{Д. Д. Реповш \footnote{Работа второго автора
поддержана Словенским иследовательским агенством, гранты BI-RU/16-18-002, P1-0292, J1-8131 и J1-7025.
}}

\vskip 0.3in

\centerline{Аннотация} В работе изучаются числовые характеристики
полиномиальных тождеств левонильпотентных алгебр. Ранее 
была предложена конструкция, позволяющая по бесконечному двоичному 
слову строить левонильпотентную ступени два алгебру с заданными 
свойствами последовательности коразмерностй. Однако класс используемых бесконечных 
слов ограничивался периодическими словами и словами Штурма. В данной работе 
предложенный ранее подход обобщается на 
значительно более общий случай. Доказано, что у любой алгебры, построенной 
по двоичному слову с субэкспоненциальной функцией комбинаторной сложности,
существует PI-экспонента и вычислено ее точное значение.

\vskip 0.5in
\section{Введение}

Пусть $F$~--- поле нулевой характеристики и $A$ --- алгебра над $F$.
С алгеброй $A$ связана целочисленная последовательность $\{c_n(A)\}$,
$n=1,2,\dots~$, определяемая ее полилинейными тождествами. В общем случае 
последовательность $\{c_n(A)\}$ может расти сверхэкспоненциально. 
Например, если $A$ --- абсолютно свободная алгебра счетного ранга, то
$\{c_n(A)\}=p(n)n!$, где $p(n)=\frac{1}{n}{2n-2\choose n-1}$ --- число  
Каталана. Если же $A$ --- свободная ассоциативная алгебра или свободная алгебра
Ли, то $\{c_n(A)\}=n!$ или $(n-1)!$ соответственно.

Тем не менее, широк класс алгебр, для которых $\{c_n(A)\}$ экспоненциально 
ограничена, т.е. $\{c_n(A)\}\le a^n$ для всех $n\ge 1$  при некоторой константе
$a$. К таким алгебрам относятся все конечномерные алгебры \cite{1},  все
ассоциативные PI-алгебры \cite{2}, бесконечномерные простые алгебры Ли
картановского типа \cite{3}, аффинные алгебры Каца-Муди \cite{z1}, \cite{z2},
алгебры Ли с нильпотентным коммутантом \cite{4} и многие другие.
В случае экспоненциальной ограниченности $\{c_n(A)\}$ последовательность корней
$\sqrt[n]{c_n(A)}$ ограничена и естественным выглядит вопрос о существовании
ее предела.

В конце 80-х годов прошлого столетия Ш. Амицур выдвинул гипотезу о том, что для
любой ассоциативной PI-алгебры $A$ предел
\begin{equation}\label{e1}
\lim_{n\to\infty}\sqrt[n]{c_n(A)}
\end{equation}
существует и является целым числом. Эта гипотеза была подтверждена в 90-х
годах \cite{5}, \cite{6}. Позднее аналогичная гипотеза нашла свое подтверждение
для конечномерных алгебр Ли \cite{7}, конечномерных йордановых алгебр \cite{8},
супералгебр Ли с нильпотентным коммутантом \cite{8a} и ряда других. Для 
некоторых бесконечномерных алгебр Ли гипотеза нашла частичное подтверждение 
---  было доказано существование предела (\ref{e1}), но он оказался дробным \cite{10a}, \cite{10b}.
существование предела (\ref{e1}) было доказано и для всех конечномерных
простых алгебр \cite{11}. При этом также были построены примеры простых
конечномерных алгебр с дробным пределом (\ref{e1}) \cite{10}. В случае существования предела
(\ref{e1}) его принято называть PI-экспонентой $A$ и обозначать как $exp(A)$.

К настоящему времени известна лишь одна работа, в которой построена серия алгебр, у которых PI-экспонента отсутствует \cite{12}. В связи с этим важную роль играет расширение класса алгебр, у которых существует PI-экспонента. 
В работе \cite{13} построено семейство неассоциативных алгебр, у которых 
PI-экспоненты пробегают все вещественные значения из бесконечного интервала
$(1;\infty)$. Все алгебры этого семейства левонильпотентны ступени два, т.е.
удовлетворяют тождеству
\begin{equation}\label{e2}
x(yz)\equiv 0.
\end{equation}
Построение примеров в \cite{13}  базировалось на комбинаторных свойствах
бесконечных двоичных слов Штурма и периодических слов. 

Недавно \cite{14}
было доказано, что любая конечнопорожденная алгебра $A$ с тождеством (\ref{e2})
имеет экспоненциально ограниченный рост последовательности коразмерностей
$\{c_n(A)\}$. Основной целью данной работы является обобщение одного из
главных результатов статьи \cite{13} --- теоремы 5.1 --- на гораздо более широкий
класс алгебр, связанных с двоичными словами. Предложено понятие {\it наклона}
для произвольного двоичного бесконечного слова $w$ и доказано, что в случае
субэкспоненциального роста комбинаторной сложности $w$ PI-экспонента соответствующей алгебры существует и явно выражается через наклон $w$.

Все основные понятия теории алгебр с тождествами читатель может найти в
\cite{15}, \cite{16}, \cite{17}. 

\section{Основные понятия и конструкции}

Обозначим через $F\{X\}$ абсолютно свободную алгебру над полем $F$ с бесконечным
множеством свободных порождающих $X$. Напомним, что полином 
$f=f(x_1,\ldots,x_n)\in F\{X\}$ называется {\it тождеством} алгебры $A$,
если $ f(a_1,\ldots,a_n)=0$ для любых $ a_1,\ldots,a_n\in A$. Совокупность
всех тождеств алгебры $A$ образует двусторонний идеал $Id(A)$ в $F\{X\}$,
устойчивый относительно всех эндоморфизмов $F\{X\}$. Обозначим через $P_n$
подпространство всех полилинейных многочленов от $ x_1,\ldots,x_n$ в
$F\{X\}$. Тогда пересечение $P_n\cap Id(A)$ состоит из всех полилинейных
тождеств $A$ степени $n$. Положим
$$
P_n(A)=\frac{P_n}{P_n\cap Id(A)},\quad c_n(A)=\dim P_n(A).
$$
Симметрическая группа $S_n$ действует естественным образом на $P_n$,
\begin{equation}\label{e3}
\sigma\circ f(x_1,\ldots,x_n)=f(x_{\sigma(1)},\ldots,x_{\sigma(n)}).
\end{equation}
Поскольку идеал $Id(A)$ устойчив относительно автоморфизма, заданного на
$F\{X\}$ формулой (\ref{e3}), то $P_n(A)$ также является $FS_n$-модулем.

Представления симметрической группы играют ключевую роль в количественной 
PI-теории, поэтому мы напомним основы понятия и конструкции, используемые
в дальнейшем. Пусть $\lambda\vdash n$ --- разбиение натурального числа $n$,
т.е. $\lambda=(\lambda_1,\ldots, \lambda_k)$, где $\lambda_1\ge\ldots
\lambda_k>0$ --- целые числа, причем $\lambda+\cdots+\lambda_k=n$.
{\it Диаграммой Юнга} $D_\lambda$ называется таблица из $n$ клеток,
Расположенных в $k$ строках. Первая строка содержит $\lambda_1$ клеток,
2-я --- $\lambda_2$ клеток, и т.д. {\it Таблицей Юнга} $T_\lambda$ называется
диаграмма $D_\lambda$ с расставленными в клетках числами от $1$  до $n$.
{\it Стабилизатором строк} $R_{T_\lambda}$ таблицы $T_\lambda$ называется
подгруппа в $S_n$, изоморфная $S_{\lambda_1}\times\cdots\times S_{\lambda_k}$
и состоящая из всех подстановок, переставляющих числа только в пределах их
строк. Аналогично, {\it стабилизатор столбцов} $C_{T_\lambda}$ состоит из
подстановок, переставляющих символы в пределах столбцов. Обозначим через
$R(T_\lambda), C(T_\lambda)$ и $e_{T_\lambda}$ следующие элементы группового
кольца $FS_n$:
$$
R(T_\lambda)=\sum_{\sigma\in R_{T_\lambda}}\sigma,\quad 
C(T_\lambda)=\sum_{\tau\in C_{T_\lambda}}({\rm sgn}~\tau)\,\tau,
\quad e_{T_\lambda}= R(T_\lambda) C(T_\lambda).
$$
Элемент $ e_{T_\lambda}$ является квазиидемпотентом, т.е.
$ e_{T_\lambda}^2=\gamma e_{T_\lambda}$, где $\gamma$ --- ненулевой
скаляр. В частности, $ C(T_\lambda) e_{T_\lambda}\ne 0$. Известно,
что левый  идеал $FS_n e_{T_\lambda}$ является минимальным и
любое неприводимое представление группы $S_n$ изоморфно ее представлению
на одном из идеалов $FS_n e_{T_\lambda}$. С основами теории представлений
симметрической группы можно познакомиться в \cite{18}, а с ее применением
в PI-теории --- в монографиях \cite{15}, \cite{16}, \cite{17}. 

При изучении  числовых характеристик, связанных с тождественными
соотношениями, принято пользоваться обозначениями из теории характеров.
Пусть $\chi_\lambda=\chi(FS_n e_{T_\lambda})$ --- характер неприводимого
модуля  $FS_n e_{T_\lambda}$. Тогда выражение
\begin{equation}\label{e4}
\chi(M)=\sum_{\lambda\vdash n} m_\lambda\chi_\lambda,
\end{equation}
где $M$ --- некоторый $FS_n$-модуль означает, что в разложении $M$ на
неприводимые компоненты
\begin{equation}\label{e5}
M=M_1\oplus\cdots\oplus M_q
\end{equation}
модуль, изоморфный $FS_n e_{T_\lambda}$, встречается $m_\lambda$ раз.
Целое неотрицательное число $m_\lambda$ называют {\it кратностью} характера
$\chi_\lambda$ в $\chi(M)$.

Пусть теперь $M$ в (\ref{e4}) --- подмодуль в $P_n$, дополнительный к
$P_n\cap Id(A)$ и изоморфный $P_n(A)$. Мы будем отождествлять  $M$ с
$P_n(A)$. Рассмотрим одно из слагаемых $M_j$ в (\ref{e5}). Оно 
порождается  полилинейным многочленом вида $ e_{T_\lambda}g$, где
$g=g(x_1,\ldots,x_n)\in P_n$. Поскольку $ e_{T_\lambda}$ --- 
квазиидемпотент, то $f= C(T_\lambda) e_{T_\lambda}g$ --- ненулевой
элемент из $M_j$. С другой стороны, если $\lambda_1^\prime,\ldots,
\lambda_t^\prime$ --- высоты столбцов диаграммы $D_\lambda$ (здесь 
$t=\lambda_1$), то множество переменных $\{x_1,\ldots,x_n\}$ разбито в
объединение непересекающихся подмножеств $X_1\cup\ldots\cup X_t$ порядков
$\lambda_1^\prime,\ldots,\lambda_t^\prime$ соответственно, таких, что $f$
кососимметричен по каждому из подмножеств $X_i, i=1,\ldots,t$. Таким
образом, мы получаем следующее утверждение.

\begin{Lemma}\label{l1}
Пусть $M_j$ --- одно из неприводимых слагаемых в (\ref{e5}) и $\chi(M_j)=
\chi_\lambda$, где $\lambda=(\lambda_1,\ldots, \lambda_k)\vdash n$. Обозначим
через $\lambda_1^\prime,\ldots,\lambda_t^\prime$ высоты столбцов $D_\lambda$.
Тогда $M_j$ порождается как $FS_n$-модуль полилинейным многочленом $f=
f(x_1,\ldots,x_n)$, таким, что $\{x_1,\ldots,x_n\}=X_1\cup\ldots\cup X_t$,
$|X_j|=\lambda_j',j=1,\ldots,t$, и $f$ кососимметричен по переменным каждого
подмножества $X_1,\ldots,X_t$. Более того, $m_\lambda\ne 0$ в (\ref{e4})
тогда и только тогда, когда $f$  не является тождеством алгебры $A$.
\end{Lemma}
\hfill $\Box$

Рассмотрим еще несколько количественных характеристик, связанных с тождествами
и представлениями $S_n$. Величина $q$ в (\ref{e5}) называется длиной модуля
$M$, а в случае $M=P_n(A)$ --- кодлиной $l_n(A)$ алгебры $A$. Ясно, что
\begin{equation}\label{e6}
l_n(A)=\sum_{\lambda\vdash n}m_\lambda,
\end{equation}
Где $m_\lambda$ --- кратности из (\ref{e4}). Если обозначить через $d_\lambda=
\deg\chi_\lambda$ размерность соответствующего представления, то
\begin{equation}\label{e7}
c_n(A)=\sum_{\lambda\vdash n}m_\lambda d_\lambda.
\end{equation}

Из (\ref{e6}) и (\ref{e7}) следует верхняя оценка
\begin{equation}\label{e8}
c_n(A)\le l_n(A) \max\{d_\lambda|~m_\lambda\ne 0\}.
\end{equation}

Вместо размерности $d_\lambda$ удобнее использовать близкую к ней
числовую характеристику. Пусть $\lambda=(\lambda_1,\ldots, \lambda_k)\vdash n$.
Обозначим
$$
\Phi(\lambda)=\frac{1}{\left( \frac{\lambda_1}{n}  \right)^{ \frac{\lambda_1}{n}}
\ldots \left( \frac{\lambda_k}{n}  \right)^{ \frac{\lambda_k}{n}}}\ .
$$
Величины $\Phi(\lambda)^n$ и $d_\lambda$ асимптотически совпадают с точностью
до полиномиального множителя. Уточним это утверждение для $k=2$. В этом
случае $\Phi(\lambda)$ является фактически функцией от одной переменной 
\begin{equation}\label{e9}
\Phi(x)=\frac{1}{x^x(1-x)^{1-x}},\quad 0<x\le\frac{1}{2}.
\end{equation}

\begin{Lemma}\label{l2}(\cite[лемма 3.3]{13})
Пусть $\lambda=(\lambda_1, \lambda_2)$ --- разбиение числа $n$. Тогда
$$
\frac{1}{\sqrt{\pi n^3}}\Phi(\beta)^n < d_\lambda<\Phi(\beta)^n,
$$
где $\beta=\frac{\lambda_2}{n}$, а $\Phi(x)$ задана формулой (\ref{e9}). 
\end{Lemma}
\hfill $\Box$

Отметим еще одно свойство функции $\Phi$.

\begin{Lemma}\label{l3}
Функция $\Phi(x)$ непрерывна в промежутке $(0;\frac{1}{2}]$ и $\Phi(a)<\Phi(b)$
для любых $0<a< b\le \frac{1}{2}$.
\end{Lemma}
\hfill $\Box$

Теперь напомним некоторые понятия комбинаторной теории бесконечных слов .
Пусть $w=w_1w_2\ldots$ --- бесконечное слово в двоичном алфавите $\{0;1\}$.
{\it Комбинаторной сложностью} слова $w$ называется функция $Comp_w(n)$,
равная числу различных подслов длины $n$, встречающихся в $w$ (см. \cite{19}).
Известно, что если $w$ --- периодическое слово, то $Comp_w(n)=const$,
начиная с некоторого $n$. В противном случае $Comp_w(n)\ge n+1$. Слово
$w$ называется {\it словом Штурма}, если $Comp_w(n)=n+1$  для всех $n$.
Наряду с медленным ростом функции сложности периодические слова и слова
Штурма обладают еще одной важной статистической характеристикой. И в том и
в другом случае существует предел
$$
\pi(w)=\lim_{n\to\infty}\frac{w_1+\cdots+w_n}{n},
$$
называемый {\it наклоном} $w$.

Чтобы обобщить результат работы \cite{13} на более широкий класс алгебр, мы
расширим понятие наклона. Пусть сначала $u=u_1\cdots u_n$ --- конечное слово 
в двоичном алфавите. Назовем {\it наклоном} $u$ величину
$$
\pi(u)=\frac{u_1+\cdots+u_n}{n},
$$
а длиной $u$ --- число $|u|=n$. Теперь для бесконечного слова $w$ положим
$$
q_n=\min\{\pi(u)|~u - ~ \hbox{подслово}~ w~ \hbox{длины} ~ n\}
$$
и определим наклон $\pi(w)$ как
$$
\pi(w)={\underline \lim}_{n\to\infty} q_n.
$$

Величина $\pi(w)$ будет задавать PI-экспоненты алгебр, соответствующих слову $w$
при гораздо более слабых условиях на функцию $Comp_w$ нежели в случае 
периодических слов и слов Штурма.

Нам понадобится следующее свойство только что определенных величин.

\begin{Lemma}\label{l3a}
Пусть $w$ --- слово с наклоном $\pi(w)=\alpha<1$. Тогда для любого натурального
$T$ в $w$ найдется подслово $u$ длины $T$ с наклоном $\pi(u)\le\alpha$.
\end{Lemma}

{\em Доказательство.} Зафиксируем произвольное достаточно малое $\varepsilon>0$. По определению
$\pi(w)$ в $w$ есть конечные подслова неограниченной длины с наклоном меньше
$\alpha+\varepsilon$. Пусть $v$ --- одно из таких слов достаточно большой
длины $m$. Разделим $m$ на $T$ с остатком,
$$
m=NT+q,
$$
и разобьем $v$ на $N+1$ подслово: $v=v_1\cdots v_{N+1}$, где последовательно
идущие подслова $v_1,\ldots,v_N$ имеют длину $T$, а $v_{N+1}$ --- подслово
длины $q$, либо пустое слово, если $q=0$. Очевидно, что при $q=0$ среди слов
$v_1,\ldots,v_N$ есть хотя бы одно с наклоном меньше $\alpha+\varepsilon$,
иначе и наклон $v$ был бы не меньше $\alpha+\varepsilon$. 

Рассмотрим теперь общий случай $0<q<T$ и предположим, что $\pi(v_i)\ge
\alpha+\varepsilon$ для всех $i=1,\ldots, N$. Наклон любого слова длины  $T$
может принимать значения только в пределах множества
$$
\{0,\frac{1}{T}, \frac{2}{T},\ldots, \frac{T-1}{T}\}.
$$
Следовательно, найдется такое целое $k$, что 
$$
\frac{k}{T}<\alpha+\varepsilon,~\frac{k+1}{T}\ge\alpha+\varepsilon.
$$
Тогда, согласно нашему предположению, в каждом из подслов $v_1,\ldots,v_N$
не менее $k+1$ единицы, и в слове $v_1\ldots v_N$ не менее $N(k+1)$ единиц.
Даже если $\pi(v_{N+1})=0$, т.е. $v_{N+1}$ состоит из одних нулей, то в $v$
не менее $N(k+1)$ единиц. Поэтому
$$
\pi(v)\ge \frac{(k+1)N}{NT+q}=\frac{k+1}{T+\frac{q}{N}}.
$$
Поскольку $\frac{k+1}{T}\ge\alpha+\varepsilon$, а $\frac{q}{N}\to 0$ при
$N\to\infty$, то условие $\pi(v)<\alpha+\varepsilon$ не выполняется при
достаточно большой длине $m$ слова $v$. Это означает, что для любого сколь угодно малого 
$\varepsilon>0$ найдется подслово $u$ длины $T$ с $\pi(u)<\alpha+\varepsilon$.
Но поскольку $\pi(u)$ для слова $u$ длины $T$ может принимать дискретный набор 
значений вида $\frac{k}{T}, 0\le k\le T-1,$ то найдется и подслово $u$ длины $T$ с
$\pi(u)\le\alpha$.

\hfill $\Box$

Лемма \ref{l3a}, в частности, означает, что верхний и нижний пределы поcледовательности 
$q_n$ совпадают, т.е. она имеет обычный предел.

\section{Алгебры двоичных слов}

Пусть $K=\{k_1,k_2,\ldots\}$ --- бесконечная целочисленная последовательность,
в которой $k_i\ge 2$ для всех $i\ge 1$. Определим алгебру $A(K)$ следующим
образом. Ее базис является бесконечным объединением
$$
\{a,b\}\cup Z_1\cup Z_2\cup\ldots,
$$
где 
$$
Z_i=\{z_j^{(i)}|1\le j\le k_i\},~i=1,2,\ldots~,
$$
а ненулевые произведения базисных элементов заданы соотношениями
\begin{equation}\label{a1}
z_2^{(i)}a= z_3^{(i)},\ldots, z_{k_i-1}^{(i)}a= z_{k_i}^{(i)},
z_{k_i}^{(i)}a= z_{1}^{(i)},~i=1,2,\ldots,
\end{equation}
$$
z_{1}^{(i)}b= z_{2}^{(i+1)},~ i=1,2,\ldots.
$$

Очевидно, что алгебра $A(K)$ удовлетворяет тождеству (\ref{e2}), а значит
ненулевыми могут быть только произведения с левонормированной расстановкой 
скобок. Поэтому мы опускаем скобки в записи одночленов. Далее, линейная
оболочка $<Z_1\cup Z_2\cup\ldots>$ является идеалом коразмерности два с
нулевым умножением. Следовательно, любой многочлен $f$ кососимметричный по
четырем переменным либо по двум наборам переменных мощности три является
тождеством алгебры $A(K)$. Отсюда в силу леммы \ref{l1} получаем следующее утверждение.

\begin{Lemma}\label{l4}
Пусть
$$
\chi_n(A(K))=\sum_{\lambda\vdash n} m_\lambda\chi_\lambda.
$$
Тогда ненулевые кратности $m_\lambda$ могут встречаться только у разбиений
$\lambda=(n), \lambda=(n-k,k), \lambda=(n-k-1,k,1)$.
\end{Lemma}
\hfill $\Box$

Теперь мы конкретизируем вид последовательности $K$. Пусть $w$ --- бесконечное 
двоичное слово и $m\ge 2$ --- целое число. Положим
\begin{equation}\label{a2}
k_i = \left\{
               \begin{array}{ll}
                   m,~ \hbox{если}~ w_i=0~, \\
                 m+1,~ \hbox{если}~ w_i=1~.
               \end{array}
          \right.
\end{equation}

Полученную алгебру, зависящую от $m$ и $w$, будем обозначать через $A(m,w)$.

\begin{Lemma}\label{l5}(\cite[лемма 4.2]{13})
Для любого $m\ge 2$ и для любого слова $w$ кодлина алгебры $A(m,w)$ при всех
$n$ удовлетворяет неравенству
$$
l_n(A(m,w)) \le 3(m+1) n^3 Comp_w(n).
$$
\end{Lemma}
\hfill $\Box$

В отличие от (\cite{13}) мы будем рассматривать алгебры не с условием 
$Comp_w(n)\le n+1$ на $w$, а с любым субэкспоненциальным ростом функции
сложности. Функцию натурального аргумента $\varphi(n)$ мы называем 
{\it субэкспоненциальной}, если
$$
\lim_{n\to\infty}\frac{\varphi(n)}{a^n}=0.
$$
Для любого вещественного $a>1$. Лемма \ref{l5} и формула (\ref{e8}) показывают,
что для слова $w$ с субэкспоненциальной функцией сложности величины
$$
\sqrt[n]{c_n(A(m,w))}\quad \hbox{и} \quad 
\sqrt[n]{\max\{d_\lambda|m_\lambda\ne 0\}}
$$ 
асимптотически близки. Заметим, что функция сложности слова $w$ --- это в
точности функция комбинаторной сложности факторного языка, состоящего из
всех конечных подслов слова $w$, а класс языков с субэкспоненциальной 
комбинаторной сложностью достаточно богат (см., например, \cite{20},
\cite{21}).

Перейдем к получению верхней оценки для коразмерностей $c_n(A(m,w))$.

\begin{Lemma}\label{l6}
Пусть $A=A(m,w)$, где $m\ge 2$, а $w$ --- бесконечное двоичное слово с
наклоном $\pi(w)=\alpha$. Пусть также $\beta=\frac{1}{m+\alpha}$. Тогда
для любого $\varepsilon>0$ существует такое $N=N(\varepsilon)$, что при всех
$n\ge N$ кратность $m_\lambda$ в разложении кохарактера $\chi_{n+1}(A)$
равна нулю, если $\frac{\lambda_2}{n}>\beta+\varepsilon$.
\end{Lemma}

{\em Доказательство.} Сначала проанализируем возможные значения полилинейных
одночленов в алгебре $A$  при подстановке базисных элементов вместо переменных.
Предположим, что найдется такой ассоциативный одночлен $g=g(a,b)$ степени
$n$, что $z_j^{(i)}g(a,b)\ne 0$ (фактически $g(a,b)$ --- это одночлен от
правых умножений на $a$ и $b$). Тогда $g$ имеет вид:
$$
g=a^{i_0}b a^{i_1}b\cdots b a^{i_{r+1}},
$$
где $0\le i_0,i_{r+1}\le m$ и $i_0,\ldots,i_r\in\{m-1,m\}$. Более того,
\begin{equation}\label{e10}
i_1+\cdots+i_r =(m-1)r+r\pi(w_{i+1}\ldots w_{i+r}),
\end{equation}
где $\pi(w_{i+1}\ldots w_{i+r})$ --- наклон подслова $w_{i+1}\ldots w_{i+r}$
слова $w$. Кроме того,
\begin{equation}\label{e11}
n=\deg g=i_0+i_{r+1}+i_1+\cdots+i_r +r+1
\end{equation}
$$
=i_0+i_{r+1}+(m-1)r+r\pi(w_{i+1}\ldots w_{i+r})+r+1,
$$
как следует из (\ref{e10}). Поскольку наклон любого слова не превосходит 1,
то из (\ref{e11}) получаем неравенство
$$
n\le 2m+(m-1)r+2r+1\le(3m+2)r.
$$
В частности,
\begin{equation}\label{e12}
r\ge\frac{n}{3m+2},
\end{equation}
т.е. $r$ растет линейно с ростом $n$. Из определения наклона $\pi(w)
=\alpha$ следует, что для любого $\delta>0$ величина
$\pi(w_{i+1}\ldots w_{i+r})$ строго больше $\alpha-\delta$ для всех
достаточно больших $n$ и $r$. Поэтому (\ref{e11}) влечет и неравенство
$$
n\ge(m-1)r+r(\alpha-\delta)+r+1 = r(m+\alpha-\delta)+1,
$$
откуда
\begin{equation}\label{e13}
\frac{\deg_bg}{\deg g}=\frac{r+1}{n}\le\frac{1}{m+\alpha-\delta}
-\frac{1}{n(m+1)}+\frac{1}{n}.
\end{equation}

Теперь рассмотрим разбиение $\lambda\vdash n+1$ с $\frac{\lambda_2}{n}
\ge\beta+\varepsilon$ и предположим, что $m_\lambda\ne 0$ в разложении
кохарактера $\chi_{n+1}(A)$. Из леммы \ref{l1}  следует, что существует полилинейный
многочлен $f=f(x_1,\ldots,x_{n+1})$, зависящий от $\lambda_2$ кососимметричных 
наборов порядка не меньше двух, не являющийся тождеством $A$. Это означает, что
существует подстановка $\varphi:\{ x_1,\ldots,x_{n+1}\}\to \{z_j^{(i)},a,b\}$,
при которой $\varphi(f)\ne 0$ и $b$ встречается среди $\varphi(x_1), \ldots, 
\varphi(x_{n+1})$ как минимум $\lambda_2-1$ раз. Но тогда $\varphi(f)$ --- 
линейная комбинация одночленов вида $z_j^{(i)}g(a,b)$, где $g$ --- ассоциативный
одночлен степени $n$, причем $\deg_bg\ge\lambda_2-1$. Следовательно,
\begin{equation}\label{e14}
\frac{\deg_bg}{\deg g} \ge\frac{1}{m+\alpha}+\varepsilon - \frac{1}{n}.
\end{equation}
Поскольку предел правой части (\ref{e13}) равен $\frac{1}{m+\alpha-\delta}$,
а предел правой части (\ref{e14}) равен $\frac{1}{m+\alpha}+\varepsilon$, 
то при фиксированном $\varepsilon$ можно подобрать такое $\delta$, что
(\ref{e13}) и (\ref{e14}) одновременно невыполнимы. Полученное противоречие завершает доказательство леммы.
\hfill $\Box$

\begin{Lemma}\label{l7}
Пусть $w$ --- бесконечное двоичное слово с наклоном $\pi(w)=\alpha$, $m\ge 2$
и $A=A(m,w)$. Тогда для любого $\nu>0$ для всех достаточно больших $n$
выполняется неравенство
$$
c_{n+1}(A)\le 3(m+1)(n+1)^4Comp_w(n+1)(\Phi(\beta)+\nu)^{n+1},
$$
где $\beta=\frac{1}{m+\alpha}$, а $\Phi(x)$ задается формулой (\ref{e9}).
\end{Lemma}

{\em Доказательство}. Зафиксируем $\nu>0$. Так как $\Phi$ непрерывна (см.
лемму \ref{l3}), то найдется такое $\varepsilon$, что $|\Phi(x)-\Phi(\beta)|
<\nu$, если $|x-\beta|<\varepsilon$. Пусть теперь $\lambda$ --- разбиение 
$n+1$ с ненулевой кратностью $m_\lambda$ в $\chi_{n+1}(A)$. По лемме
\ref{l6} имеем: $\frac{\lambda_2}{n}\le\beta+\varepsilon$. Обозначим
$\rho=\frac{\lambda_2}{n}$. По лемме \ref{l4} у разбиения $\lambda$
третья компонента $\lambda_3$ равна 0 или 1. Если $\lambda=(\lambda_1,
\lambda_2)$, т.е. $\lambda_3=0$, то
$$
d_\lambda=\deg\chi_\lambda \le \Phi(\rho)^{n+1} 
\le(\Phi(\beta)+\nu)^{n+1}
$$
по лемме \ref{l2}. Для одноэлементного разбиения $\lambda=(n+1)$ аналогичное 
неравенство очевидно. Пусть $\lambda_3=1$. Из формулы крюков для размерностей 
неприводимых представлений группы $S_n$ (см. \cite{18}) следует, что
$$
d_\lambda=d_\mu \frac{\lambda_2(\lambda_1+1)(n+1)}{(\lambda_2+1)(\lambda_1+2)}
<(n+1)d_\mu,
$$
где $\mu=(\lambda_1, \lambda_2)\vdash n$. Следовательно, $d_\lambda <(n+1)
\Phi(\rho)^n$ по лемме \ref{l2}. Таким образом, неравенство
\begin{equation}\label{e14a}
d_\lambda <(n+1)(\Phi(\beta)+\nu)^{n+1}
\end{equation}
выполняется для всех $\lambda\vdash(n+1)$ с $m_\lambda \ne 0$. Из (\ref{e14a}),
(\ref{e8}) и леммы \ref{l5} мы получаем
$$
c_{n+1}(A)\le l_{n+1}(A)(n+1)(\Phi(\beta)+\nu)^{n+1}\le
$$
$$
3(m+1)(n+1)^4Comp_w(n+1)(\Phi(\beta)+\nu)^{n+1},
$$
и лемма доказана.
\hfill $\Box$

\section{Основной результат}

Перейдем к получению нижней оценки роста коразмерностей алгебры $A(m,w)$.

\begin{Lemma}\label{l8}
Пусть $w$ --- слово с наклоном $\pi(w)=\alpha$, $m\ge 2$ и $A=A(m,w)$. Тогда для 
любого $0<\varepsilon<\frac{1}{m+\alpha}$ существуют $r_0$ и 
возрастающая последовательность $n_{r_0}, n_{r_0+1},\ldots~$, такие, что
\begin{equation}\label{e14aa}
c_{n_r+1}(A)\ge\frac{1}{\sqrt{\pi n_r^3}}
\Phi(\frac{1}{m+\alpha}-\varepsilon)^{n_r},\quad r\ge r_0.
\end{equation}
Кроме того, $n_{r+1}-n_r\le 2m+1$ для всех $r\ge r_0$.
\end{Lemma}

{\em Доказательство}. Выберем $r_0$ таким что $\varepsilon r_0>m$. Пусть теперь
$r\ge r_0$ --- любое целое число. По лемме \ref{l3a} в $w$ найдется подслово
$v=w_{i+1}\ldots w_{i+r}$ с наклоном $\pi(v)\le\alpha$. Как и в лемме \ref{l6}
в алгебре $A$ найдется ненулевое произведение вида $z_j^{(i)}g(a,b)$, где
$g=g(a,b)$ --- ассоциативный одночлен от правых умножений на $a$ и $b$,
$$
g=a^{i_0}ba^{i_1}b\cdots ba^{i_r},
$$
причем $0\le i_0\le m$, а $i_1,\ldots,i_r\in\{m-1,m\}$. Кроме того,
$$
i_1+\cdots+i_r=(m-1+\pi(v))r
$$
и
$$
n=\deg g=i_0+i_1+\cdots+i_r+r=i_0+r(m+\pi(v))\le m+r(m+\alpha),
$$
а степень $g$ по $b$ равна $r$. Из предыдущего неравенства и выбора $r$ следует, что 
\begin{equation}\label{e15}
\frac{r}{n}\ge\frac{1}{m+\alpha}-\varepsilon.
\end{equation}

Рассмотрим в свободной алгебре $F\{X\}$ левонормированный одночлен
$$
f=zx_1^0\cdots x_{i_0}^0 y_1 x_1^1\cdots x_{i_1}^1 y_2\cdots y_r
x_1^r\cdots x_{i_r}^r.
$$
При подстановке $\varphi: X\to A$, $\varphi(z)=z_j^{(i)}, \varphi(x_q^p)=a,
\varphi(y_s)=b$ для всех допустимых $p,q,s$ значение $\varphi(f)$ равно
$z_j^{(i)}g(a,b)$, т.е. $\varphi(f)\ne 0$.  Теперь проальтернируем $f$ по парам
$\{x^k_1,y_k\}, 1\le k\le r$, и получим многочлен $\widetilde f$. Таблица 
умножения базисных элементов алгебры $A$ показывает, что $\varphi(\widetilde f)
=\varphi(f)\ne 0$, т.е. $\widetilde f$ не является тождеством алгебры $A$.

Рассмотрим действие симметрической группы $S_n$ на переменные $\{x_p^q,y_s\}$ и
$FS_n$-подмодуль в $P_{n+1}$, порожденный многочленом $\widetilde f$. Структура
квазиидемпотентов кольца $FS_n$ показывает, что в разложении $FS_n\widetilde f$
на неприводимые компоненты могут возникнуть только подмодули с характером
$\chi_\lambda$, где $\lambda=(n-t,t)$ и $t\ge r$. Но тогда для одного из таких
разбиений $\lambda$ получаем
$$
c_{n+1}(A) \ge d_\lambda\ge\frac{1}{\sqrt{\pi n^3}}\Phi(\frac{t}{n})^n
$$
в силу леммы \ref{l2}. Так как $\frac{t}{n}\ge\frac{r}{n}$, то
$$
c_{n+1}(A) \ge\frac{1}{\sqrt{\pi n^3}}\Phi(\frac{1}{m+\alpha}-\varepsilon)^n
$$
с учетом (\ref{e15}) и леммы \ref{l3}.

Именно это значение $n=i_0+i_1+\cdots+i_r+r=i_0+r(m+\pi(w_{i+1}\ldots w_{i+r}))$
мы и возьмем в качестве $n_r$ для выбранного $r\ge r_0$.

Оценим разность $n_{r+1}-n_r$. Заметим сначала, что $n_{r+1}>n_r$ в силу
выбора слов $w_{t+1}\ldots w_{t+r+1}$ и $w_{i+1}\ldots w_{i+r}$ с минимальным
наклоном. Величина $n_r$ задается подсловом $w_{i+1}\ldots w_{i+r}$ в 
$w$ длины $r$ с минимальным наклоном, а произведение $r\pi(w_{i+1}\ldots w_{i+r})$ равно количеству единиц среди $w_{i+1},\ldots,
w_{i+r}$. Аналогично, для $n_{r+1}$ есть подслово $w_{t+1}\ldots w_{t+r+1}$ с
наименьшим  количеством единиц,
$$
n_{r+1}=i_0^\prime+(r+1)(m+\pi(w_{t+1}\ldots w_{t+r+1}))
$$
и поэтому
$$
n_{r+1}-n_r\le 2m+(r+1)\pi(w_{t+1}\ldots w_{t+r+1})
-r\pi(w_{i+1}\ldots w_{i+r}).
$$
Но в силу выбора подслов $w_{t+1}\ldots w_{t+r+1}$ и $w_{i+1}\ldots w_{i+r}$
количество единиц в них либо совпадает, либо отличается на 1. Поэтому
$n_{r+1}-n_r\le 2m+1$, и лемма доказана.

\hfill $\Box$

Докажем теперь основной результат работы, обобщающий теорему 5.1 из \cite{13}.

\begin{Theorem}\label{t1}
Пусть $w$ --- бесконечное двоичное слово с наклоном $\pi(w)=\alpha$ и с
функцией сложности $Comp_w(n)$ субэкспоненциального роста, и $m\ge 2$
--- целое число. Тогда у алгебры $A=A(m,w)$, заданной соотношениями (\ref{a1}),
(\ref{a2}), существует PI-экспонента, причем
\begin{equation}\label{e15a}
exp(A)= \Phi(\frac{1}{m+\alpha}),
\end{equation}
где
$$
\Phi(x)=\frac{1}{x^x(1-x)^{1-x}}.
$$
\end{Theorem}
{\em Доказательство.} Обозначим для удобства $a_n=\sqrt[n]{c_n(A)}$ и покажем,
что верхний и нижний пределы последовательности $a_n$ совпадают и равны
$\Phi(\frac{1}{m+\alpha})$.

Сначала оценим нижний предел. Для этого докажем, что для любого сколь угодно малого $\varepsilon>0$ при всех достаточно больших $n$ выполняется неравенство
\begin{equation}\label{e16}
c_{n+1}(A)\ge \frac{1}{2^{2m+1}\sqrt{\pi n^3}}
\Phi(\frac{1}{m+\alpha}-\varepsilon)^n.
\end{equation}

По лемме \ref{l8} для заданного $\varepsilon$ существует последовательность 
индексов $n_r, r=r_0, r_0+1,\ldots~,$ для которых выполняется неравенство 
(\ref{e14aa}). Для любого $n\ge n_{r_0}$ найдется такое $r$, что 
$n_r\le n< n_{r+1}$ и $n-n_r< 2m+1$.

Заметим, что правый аннулятор любого элемента $x\ne 0$ из $A$ равен нулю.
Отсюда немедленно следует, что $c_{t+1}(A) \ge c_t(A)$ для любого $t$.
Поэтому
\begin{equation}\label{e17}
c_{n+1}(A)\ge c_{n_r+1}(A)\ge \frac{1}{\sqrt{\pi n^3_{n_r}}}
\Phi(\frac{1}{m+\alpha}-\varepsilon)^{n_{r}}.
\end{equation}
Так как $n\ge n_r$, $\Phi(x)\le 2$ и $n-n_r\le 2m+1$, то правая часть
(\ref{e17}) не меньше чем
$$
\frac{1}{2^{2m+1}\sqrt{\pi n^3}}\Phi(\frac{1}{m+\alpha}-\varepsilon)^n,
$$
что и доказывает соотношение (\ref{e16}).

Поскольку
$$
\left( 2^{2m+1}\sqrt{\pi n^3}\right)^{-\frac{1}{n}}\to 1
\quad \hbox{при}\quad n\to\infty,
$$
то соотношение (\ref{e16}) влечет неравенство
$$
\underline{\lim}_{n\to\infty} a_n \ge \Phi(\frac{1}{m+\alpha}-\varepsilon)
\quad\hbox{для любого}\quad\varepsilon>0,
$$
следовательно,
\begin{equation}\label{e18}
\underline{\lim}_{n\to\infty} a_n \ge \Phi(\frac{1}{m+\alpha}).
\end{equation}

Для получения оценки верхнего предела последовательности $\{a_n\}$ воспользуемся 
леммой \ref{l7}. Поскольку $Comp_w(n)$ --- функция субэкспоненциального роста,
предел величины
$$
\left(
3(m+1)(n+1)^4Comp_w(n+1)
\right)^{\frac{1}{n}}
$$
при $n\to\infty$ равен единице. Следовательно, по лемме \ref{l7}
\begin{equation}\label{e19}
\overline{\lim}_{n\to\infty} a_n \le\Phi(\frac{1}{m+\alpha}).
\end{equation}

Неравенства (\ref{e18}), (\ref{e19}) означают существование обычного предела
последовательности $\{a_n\}$, т.е. PI-экспоненты алгебры $A(m,w)$, а также 
равенства (\ref{e15a}).

\hfill $\Box$

В заключение выскажем гипотезу о существовании PI-экспоненты у любой 
левонильпотентной ступени два алгебры с конечным числом порождающих. Важными 
частными случаями являются конечномерные алгебры и относительно свободные 
алгебры конечного ранга. Для них поставленный вопрос также открыт. Условие 
конечной порожденности  существенно, как показывает контрпример 
\cite{14}.